\documentclass[aip,
amsmath,amssymb,
preprint,%
numerical,
]{revtex4-1}
\usepackage{graphicx}
\usepackage{dcolumn}
\usepackage{bm}
\usepackage{hyperref}
\usepackage[T1]{fontenc}
\usepackage[utf8]{inputenc}
\usepackage[utf8]{inputenc}
\usepackage[T1]{fontenc}
\usepackage{mathptmx}
\usepackage{etoolbox}
\newcommand{\RNum}[1]{\uppercase\expandafter{\romannumeral #1\relax}}
\usepackage{array,xcolor,colortbl}
\usepackage{multirow}
\makeatletter
\def\@email#1#2{%
	\endgroup
	\patchcmd{\titleblock@produce}
	{\frontmatter@RRAPformat}
	{\frontmatter@RRAPformat{\produce@RRAP{*#1\href{mailto:#2}{#2}}}\frontmatter@RRAPformat}
	{}{}
}%
\makeatother
\begin{document}

	\title{Effect of rotation on anisotropic scattering suspension of phototactic algae}
	\author{S. K. Rajput}
	\email{shubh.iiitj@gmail.com}
	\affiliation{ 
		Department of Mathematics, PDPM Indian Institute of Information Technology Design and Manufacturing,
		Jabalpur 482005, India.
	}%
	
	

	\begin{abstract}
	In this article, the effect of rotation on the onset of phototactic bioconvection is investigated using linear stability theory for a suspension of forward-scattering phototactic algae in this article. The suspension is uniformly illuminated by collimated flux. The bio-convective instability is characterized by an unstable mode of disturbance that transitions from a stationary (overstable) to an overstable (stationary) state as the Taylor number varies under fixed parameters. It is also observed that the suspension has significant stabilizing effect due to rotation of the system.
	\end{abstract}
	
	\maketitle

	\section{INTRODUCTION}
	Bioconvection is a widespread phenomenon in fluid dynamics that involves the formation of patterns in suspensions of living microorganisms, such as algae and bacteria~\cite{20platt1961,21pedley1992,22hill2005,23bees2020,24javadi2020}. The term "bioconvection" was first introduced by Platt in 1961. Typically, these microorganisms have a higher density than the surrounding fluid on a small scale and collectively swim in an upward direction. However, there are cases in natural environments where density differences are not necessary for pattern formation, and the microorganisms exhibit different behaviours. Non-living microorganisms do not exhibit pattern formation behaviour. Microorganisms respond to external stimuli and exhibit behavioural changes in their swimming direction, known as taxis. Examples of taxis include gravitaxis, which is influenced by gravitational acceleration, gyrotaxis, which is influenced by both gravitational acceleration and viscous torque when microorganisms have bottom-heavy structures, and phototaxis, which is the movement of microorganisms towards or away from a light source. This article focuses specifically on the effect of phototaxis.   

The phenomenon of bioconvection patterns in algal suspensions has been the subject of extensive research, highlighting the significant role of light in shaping fluid dynamics and cellular distribution~\cite{1wager1911,2kitsunezaki2007}. When exposed to intense light in well-mixed cultures of photosynthetic microorganisms, dynamic patterns can either persist or undergo disruption. Strong light can hinder the formation of stable patterns, and it can also disrupt pre-existing ~\cite{3kessler1985,4williams2011,5kessler1989} The response of algae to light of varying intensities is a key determinant in the modulation of bioconvection patterns. Moreover, the interaction of light with microorganisms involves processes of light absorption and scattering. Algal light scattering can be characterized as either isotropic or anisotropic, the latter further classified into forward and backward scattering based on cell properties such as size, shape, and refractive index. Algae, due to their size, predominantly exhibit forward scattering of light within the visible wavelength range.

The phototaxis model proposed by Rajput is being employed to investigate the bioconvection system. In this study, the suspension of phototactic microorganisms is assumed to undergo rotation around the z-axis at a constant angular velocity, while being illuminated from above by collimated flux. Given that many motile algae rely on photosynthesis for their energy needs, they exhibit distinct phototactic behavior. To ensure a comprehensive understanding of their response to light, it is essential to examine the phototaxis model in the context of a rotating medium. Unlike previous models that considered scattering effects in an isotropic manner, our approach incorporates anisotropic scattering. This allows for a more precise representation of the swimming behaviour of microorganisms as they react to light stimuli. By considering the influence of anisotropic scattering, we aim to enhance the accuracy of the model and provide valuable insights into the intricate dynamics of the bioconvection phenomenon.

\begin{figure}[!ht]
	\centering
	\includegraphics[width=14cm]{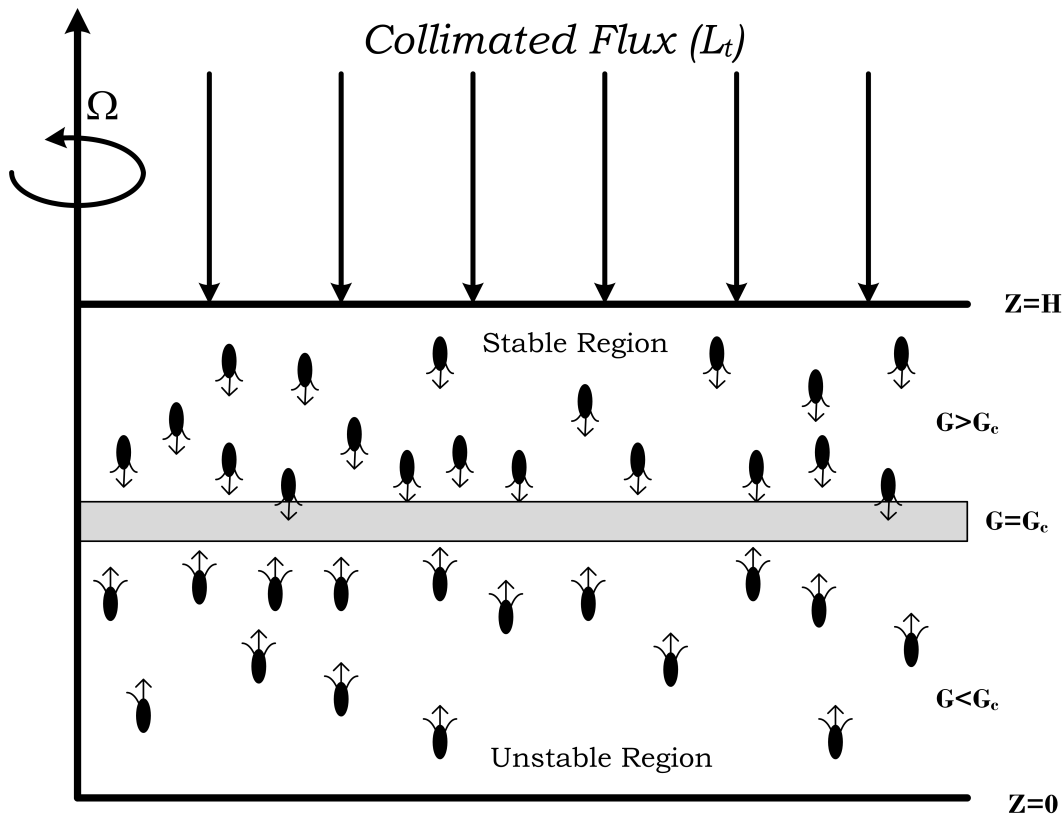}
	\caption{\footnotesize{Formation of sublayer in rotating medium.}}
	\label{fig1}
\end{figure}
In this experimental setup, the mathematical formulation of the problem involves considering the behavior of a forward scattering algal suspension confined between two parallel boundaries in the (y, z) plane. The suspension is assumed to be uniformly illuminated from above by collimated flux. The goal is to study the behavior of the suspension under the influence of light.

To analyze the system, a realistic phototaxis model is employed, which takes into account the forward scattering of light by the algae. The equilibrium solution of the system is determined by considering a steady-state where the flow velocity is zero, indicating no fluid motion. In this state, the balance between phototaxis and diffusion leads to the formation of a horizontal sublayer where cells accumulate, dividing the suspension into two distinct regions.

The critical intensity $G_c$ plays a crucial role in determining the position of the sublayer. Above the sublayer ($G > G_c$), the light intensity is sufficient to suppress cell movement, while below the sublayer ($G < G_c$), cells exhibit upward motion towards the sublayer due to phototaxis.

The lower region beneath the sublayer is considered unstable and has the potential to exhibit fluid motion and penetrate the stable upper region through penetrative convection if the system becomes unstable. Penetrative convection is a phenomenon observed in various convection problems.

Over the years, researchers have made significant advancements in modeling phototaxis and bioconvection, employing various approaches such as linear stability theory and numerical simulations. For instance, Vincent and Hill~\cite{12vincent1996} investigated the initiation of phototactic bioconvection and identified non-stationary and non-oscillatory disturbance modes. Ghorai and Hill~\cite{10ghorai2005} conducted a two-dimensional simulation of bioconvective flow patterns based on Vincent and Hill's model, albeit without considering algae scattering. Ghorai et al.~\cite{7ghorai2010} examined the onset of bioconvection in an isotropic scattering suspension and observed a steady-state profile with a sublayer located at different depths due to isotropic scattering. Panda and Ghorai~\cite{14panda2013} conducted a linear stability analysis in a forward scattering algal suspension, observing a transition from non-oscillatory to non-stationary modes during bioconvective instability caused by forward scattering, while Panda and Singh~\cite{11panda2016} numerically explored the linear stability of PBC using Vincent and Hill's continuum model in a two-dimensional setting. Previous studies did not account for the effects of diffuse irradiation. Panda et al.~\cite{15panda2016} proposed a model that examined how diffuse radiation influenced isotropic scattering algal suspensions, revealing that diffuse irradiation significantly stabilizes such suspensions. Panda~\cite{8panda2020} further investigated the impact of forward anisotropic scattering on the onset of PBC, considering both diffuse and collimated irradiation. However, these studies did not account for the influence of oblique collimated flux. Panda et al.~\cite{16panda2022} introduced oblique collimated flux and discovered that bioconvective solutions switch from non-oscillatory to overstable states and vice versa during bioconvective instability induced by oblique collimated flux.  Kumar~\cite{17kumar2022} examined the effect of oblique collimated irradiation on isotropic scattering algal suspensions, revealing the existence of both stationary and overstable solutions within certain parameter ranges. Kumar~\cite{40kumar2023} investigated the effect of rotation on the phototactic bioconvection in a non-scattering medium and found a significant stabilizing impact due to rotation. More recently, Panda and Rajput~\cite{41rajput2023} investigate the combined impact of oblique and diffuse irradiation on the onset of bioconvection. 

However, there is still a need to investigate phototactic bioconvection in a rotating forward scattering algal suspension using a realistic phototaxis model. Therefore, the objective of this study is to examine bioconvective instability in such a system and provide insights into the complex dynamics of bioconvection in a light-induced pattern formation in a dilute algal suspension.

The article follows a structured approach, starting with the mathematical formulation of the problem. It derives the equilibrium solution and perturbs the governing system of bioconvection with small disturbances. The linear stability problem is formulated and solved using numerical methods. The obtained results are then presented and discussed. The article concludes by highlighting the novelty and significance of the proposed model in advancing the understanding of bioconvective instability in rotating forward scattering algal suspensions.

	\section{MATHEMATICAL FORMULATION}
The research focuses on studying the behaviour of a forward scattering algal suspension within the (y, z) plane, which is confined between two parallel boundaries. The boundaries are assumed to be infinite in extent, and there is no reflection of light from the top and bottom boundaries. This setup allows for the investigation of the suspension's behaviour under the influence of light.

The suspension is uniformly illuminated from above by collimated flux. Collimated flux refers to light that is travelling in parallel rays, indicating that the incident light rays are not scattered before reaching the suspension. This uniform illumination from above serves as the driving force for the system and influences the phototactic behaviour of the algae.

	\subsection{THE AVERAGE SWIMMING DIRECTION}
	
	The Radiative Transfer Equation (RTE) is utilized to determine the light intensity in the system. This equation provides a mathematical description of how light propagates and interacts with the medium. The RTE can be expressed as follows 
	\begin{equation}\label{1}
		\frac{dL(\boldsymbol{x},\boldsymbol{r})}{dr}+(a+\sigma_s)L(\boldsymbol{x},\boldsymbol{r})=\frac{\sigma_s}{4\pi}\int_{0}^{4\pi}L(\boldsymbol{x},\boldsymbol{r'})\Xi(\boldsymbol{r},\boldsymbol{r'})d\Omega',
	\end{equation}
	where $a$, and $\sigma_s$ are the absorption coefficient, and scattering coefficients respectively. $\Xi(\boldsymbol{r},\boldsymbol{r'})$ is the scattering phase function, which describes how light is distributed in different directions when scattered from $\boldsymbol{r'}$ to $\boldsymbol{r}$. In this study, we assume that the scattering phase function follows linear anisotropy with azimuthal symmetry, specifically $\Xi(\boldsymbol{r},\boldsymbol{r'})=1+A\cos{\theta}\cos{\theta'}$. Here, $A$ represents the anisotropic coefficient, which determines the degree of anisotropy. Positive values of $A$ indicate forward scattering, negative values indicate backward scattering and $A=0$ represents isotropic scattering.

In the given model, the light intensity on the top of the suspension is described by:

\begin{equation*}
	L(\boldsymbol{x}_b,\boldsymbol{s}) = L_t\delta(\cos{\theta}-\cos{\theta_0}),
\end{equation*}

where $\boldsymbol{x}_b=(x,y,H)$ represents the location on the top boundary surface. $L_t$ is the magnitude of the collimated flux, $\theta$ is the zenith angle, and $\theta_0$ represents a specific value of the zenith angle.

Introducing two variables, $a=\alpha n(\boldsymbol{x})$ and $\sigma_s=\beta n(\boldsymbol{x})$, in terms of the scattering albedo $\omega=\sigma_s/(a+\sigma_s)$, the Radiative Transfer Equation (RTE) can be expressed as:

\begin{equation}\label{2}
\frac{dL(\boldsymbol{x},\boldsymbol{r})}{dr}+(\alpha+\beta)nL(\boldsymbol{x},\boldsymbol{r})=\frac{\beta n}{4\pi}\int_{0}^{4\pi}L(\boldsymbol{x},\boldsymbol{r'})(A\cos{\theta}\cos{\theta'})d\Omega',
\end{equation}

where $L(\boldsymbol{x},\boldsymbol{r})$ represents the light intensity at position $\boldsymbol{x}$ along the ray direction $\boldsymbol{r}$, $\alpha$ is the absorption coefficient, $\beta$ is the scattering coefficient, $n$ is the cell concentration, and $\Omega'$ represents the solid angle.

The total intensity at a fixed point $\boldsymbol{x}$ in the medium is given by:

\begin{equation*}
	G(\boldsymbol{x}) = \int_0^{4\pi}L(\boldsymbol{x},\boldsymbol{r})d\Omega,
\end{equation*}

and the radiative heat flux is given by:

\begin{equation}\label{3}
\boldsymbol{q}(\boldsymbol{x}) = \int_0^{4\pi}L(\boldsymbol{x},\boldsymbol{r})\boldsymbol{r}d\Omega.
\end{equation}

The average swimming velocity of a microorganism is defined as:

\begin{equation*}
	\boldsymbol{W}_c = W_c<\boldsymbol{p}>.
\end{equation*}

The mean swimming orientation of the cell, denoted as $<\boldsymbol{p}>$, is determined by employing the following calculation. In this context, $W_c$ denotes the average cell swimming speed.

\begin{equation}\label{4}
<\boldsymbol{p}> = -M(G)\frac{\boldsymbol{q}}{|\boldsymbol{q}|},
\end{equation}

where $M(G)$ is the taxis function that characterizes the response of algae cells to light. It takes the form:

\begin{equation*}
	M(G)=\left\{\begin{array}{ll}\geq 0, & \mbox{ } G(\boldsymbol{x})\leq G_{c}, \\
		< 0, & \mbox{ }G(\boldsymbol{x})>G_{c}.  \end{array}\right. 
\end{equation*}
	
The taxis function determines the behaviour of the cells in response to the light intensity.	
	
	\subsection{GOVERNING EQUATIONS WITH BOUNDARY CONDITIONS}
	
In the proposed model of the suspension of phototactic microorganisms, several equations govern the system. These equations, along with the corresponding boundary conditions, are as follows:

Continuity equation:
\begin{equation}\label{5}
\boldsymbol{\nabla}\cdot \boldsymbol{U}=0,
\end{equation}

Momentum equation under the Boussinesq approximation:
\begin{equation}\label{6}
\rho\left(\frac{\partial \boldsymbol{U}}{\partial t}+(\boldsymbol{U}\cdot\nabla )\boldsymbol{U}+2\boldsymbol{\Omega}\times \boldsymbol{U}\right)=-\boldsymbol{\nabla} P+\mu\nabla^2\boldsymbol{U}-nVg\Delta\rho\hat{\boldsymbol{z}},
\end{equation}

Cell conservation equation:
\begin{equation}\label{7}
\frac{\partial n}{\partial t}=-\boldsymbol{\nabla}\cdot \boldsymbol{B},
\end{equation}

Total cell flux equation:
\begin{equation}\label{8}
\boldsymbol{B}=nU+nW_c<\boldsymbol{p}>-\boldsymbol{D}\boldsymbol{\nabla} n,
\end{equation}

where $U$ is the average fluid velocity, $n$ is the number of algal cells per unit volume, $V$ is the volume of the cells, $\rho$ is the density of water, $\Delta\rho/\rho$ represents the small variation in density of the cells, $\boldsymbol{\Omega}$ is the angular velocity, $\mu$ is the dynamic viscosity of the suspension, $\boldsymbol{D}$ is the cell diffusivity, $DI$ represents the constant diffusivity assumption, $g$ is the acceleration due to gravity, and $\boldsymbol{p}$ is the direction of cell swimming.

The boundary conditions for the system are as follows:

For the lower and upper boundaries:
\begin{equation}\label{9}
\boldsymbol{U}\cdot\hat{\boldsymbol{z}}=0\quad \text{at}\quad z=0,H,
\end{equation}
\begin{equation}\label{10}
\boldsymbol{B}\cdot\hat{\boldsymbol{z}}=0\quad \text{at}\quad z=0,H.
\end{equation}

For the rigid boundary:
\begin{equation}\label{11}
\boldsymbol{U}\times\hat{\boldsymbol{z}}=0\quad \text{at}\quad z=0,
\end{equation}

For the stress-free boundary:
\begin{equation}\label{12}
\frac{\partial^2}{\partial z^2}(\boldsymbol{U}\cdot\hat{\boldsymbol{z}})=0\quad \text{at}\quad z=H.
\end{equation}

Additionally, the upper boundary is subjected to collimated flux, which leads to specific conditions for the intensity at the boundaries:

For the upper boundary:
\begin{subequations}
	\begin{equation}\label{13a}
	L(x, y, z = 0, \theta, \phi)=L_t\delta(\boldsymbol{r}-\boldsymbol{r_0}),\quad \left(\frac{\pi}{2}\leq\theta\leq\pi\right),
	\end{equation}
	\begin{equation}\label{13b}
	L(x, y, z = 0, \theta, \phi) =0,\quad \left(0\leq\theta\leq\frac{\pi}{2}\right).
	\end{equation}
\end{subequations}

In these equations, $L_t$ represents the intensity of light at a specific location $\boldsymbol{r_0}$, and $\delta(\boldsymbol{r}-\boldsymbol{r_0})$ denotes the Dirac delta function. These boundary conditions specify the behaviour of the fluid velocity, cell concentration, and light intensity at the boundaries of the system.

To express the governing equations in a dimensionless form, we introduce the following scales: Length $H$, time $H^2/D$, velocity scale $D/H$, pressure scale $\mu D/H^2$, concentration $\bar{n}$. Using these scales, we can non-dimensionalize the governing equations 
		\begin{equation}\label{14}
		\boldsymbol{\nabla}\cdot\boldsymbol{U}=0,
	\end{equation}
	\begin{equation}\label{15}
		S_c^{-1}\left(\frac{\partial \boldsymbol{U}}{\partial t}+(\boldsymbol{U}\cdot\nabla )\boldsymbol{U}\right)+\sqrt{T_a}(\hat{z}\times\boldsymbol{U})=-\nabla P_{e}-Rn\hat{\boldsymbol{z}}+\nabla^{2}\boldsymbol{U},
	\end{equation}
	\begin{equation}\label{16}
		\frac{\partial{n}}{\partial{t}}=-\boldsymbol{\nabla}\cdot\boldsymbol{B}=-{\boldsymbol{\nabla}}\cdot[\boldsymbol{n{\boldsymbol{U}}+nV_{c}<{\boldsymbol{p}}>-{\boldsymbol{\nabla}}n.}]
	\end{equation}
	
	In the above equations, $S_c^{-1}=\nu/D$ represents the Schmidt number, $V_c$ denotes the dimensionless swimming speed as $V_c=W_cH/D$, $R=\tilde{n}V g\Delta{\rho}H^{3}/\nu\rho{D}$ is the Rayleigh number, and $T_a=4\Omega^2H^4/\nu^2$ is the Taylor number.

After non-dimensionalization, the boundary conditions can be expressed in a dimensionless form as follows:

\begin{equation}\label{17}
U\cdot\hat{z}=0 \quad \text{at} \quad z=0,1,
\end{equation}
\begin{equation}\label{18}
B\cdot\hat{z}=0 \quad \text{at} \quad z=0,1,
\end{equation}
\begin{equation}\label{19}
U\times\hat{z}=0 \quad \text{at} \quad z=0,
\end{equation}
\begin{equation}\label{20}
\frac{\partial^2}{\partial z^2}(U\cdot\hat{z})=0 \quad \text{at} \quad z=1.
\end{equation}	
	
	The nondimensional Radiative Transfer Equation (RTE) can be expressed as:
	
	\begin{equation}\label{21}
	\frac{dL(\boldsymbol{x},\boldsymbol{r})}{dr}+\kappa nL(\boldsymbol{x},\boldsymbol{r})=\frac{\sigma n}{4\pi}\int_{0}^{4\pi}L(\boldsymbol{x},\boldsymbol{r'})(A\cos{\theta}\cos{\theta'})d\Omega',
	\end{equation}
	
	where $\kappa=(\alpha+\beta)\bar{n}H$ and $\sigma=\beta\bar{n}H$ are the nondimensional extinction coefficient and scattering coefficient, respectively. The scattering albedo $\omega=\sigma/(\kappa+\sigma)$ represents the scattering efficiency of microorganisms. The RTE describes the balance between absorption and scattering of light in the medium.
	
	Alternatively, the RTE can be formulated using direction cosines, which provide a convenient way to describe the direction of light propagation:
	
	\begin{equation}\label{23}
	\xi\frac{dL}{dx}+\eta\frac{dL}{dy}+\nu\frac{dL}{dz}+\kappa nL(\boldsymbol{x},\boldsymbol{r})=\frac{\omega\kappa n}{4\pi}\int_{0}^{4\pi}L(\boldsymbol{x},\boldsymbol{r'})(A\cos{\theta}\cos{\theta'})d\Omega',
	\end{equation}
	
	In the dimensionless form, the intensity at the boundaries is given by:
	
	\begin{subequations}
		\begin{equation}\label{24a}
		L(x, y, z = 1, \theta, \phi)=L_t\delta(\boldsymbol{r}-\boldsymbol{r_0}),\qquad \left(\frac{\pi}{2}\leq\theta\leq\pi\right),
		\end{equation}
		\begin{equation}\label{24b}
		L(x, y, z = 0, \theta, \phi) =0,\qquad \left(0\leq\theta\leq\frac{\pi}{2}\right).
		\end{equation}
	\end{subequations}

	\section{THE STEADY STATE SOLUTION}
	
The equilibrium state of the system is characterized by the following conditions:

\begin{equation*}
	\boldsymbol{U}=0,~~~n=n_s(z),~~~\zeta_s=\nabla\times U=0,~~~\text{and}~~~L=L_s(z,\theta).
\end{equation*}

At equilibrium, the total intensity $G_s$ and radiative flux $\boldsymbol{q}_s$ can be expressed as:

\begin{equation*}
	G_s=\int_0^{4\pi}L_s(z,\theta)d\Omega, \quad \text{and} \quad \boldsymbol{q}_s=\int_0^{4\pi}L_s(z,\theta)\boldsymbol{s}d\Omega.
\end{equation*}

The $\boldsymbol{q}_s$ vector has zero components in $x$ and $y$ directions due to the independence of $L_s^d(z,\theta)$ on $\phi$. Therefore, $\boldsymbol{q}_s$ can be written as $\boldsymbol{q}_s=-q_s\hat{\boldsymbol{z}}$, where $q_s$ represents the magnitude of $\boldsymbol{q}_s$.

The equation governing $L_s$ at equilibrium is given by:

\begin{equation}\label{26}
\frac{dL_s}{dz}+\frac{\kappa n_sL_s}{\nu}=\frac{\omega\kappa n_s}{4\pi\nu}(G_s(z)-Aq_s\nu).
\end{equation}

The equilibrium intensity $L_s$ can be decomposed into two components: the collimated part $L_s^c$ and the diffuse part caused by scattering $L_s^d$. The equation governing the collimated component $L_s^c$ is:

\begin{equation}\label{27}
\frac{dL_s^c}{dz}+\frac{\kappa n_sL_s^c}{\nu}=0,
\end{equation}

with the boundary condition:

\begin{equation}\label{28}
L_s^c(1, \theta) =L_t\delta(\boldsymbol{r}-\boldsymbol{r}_0), \quad \theta\in[\pi/2,\pi],
\end{equation}

where $L_t$ represents the total intensity at the point source location $\boldsymbol{r}_0$.

Solving the governing equation for $L_s^c$ with the given boundary condition yields the expression for $L_s^c$:

\begin{equation}\label{29}
L_s^c=L_t\exp\left(\int_z^1\frac{\kappa n_s(z')}{\nu}dz'\right)\delta(\boldsymbol{r}-\boldsymbol{r}_0).
\end{equation}

The diffuse part $L_s^d$ is governed by the equation:

\begin{equation}\label{30}
\frac{dL_s^d}{dz}+\frac{\kappa n_sL_s^d}{\nu}=\frac{\omega\kappa n_s}{4\pi\nu}(G_s(z)-Aq_s\nu),
\end{equation}

with the boundary conditions:

\begin{equation}\label{31a}
L_s^d(1, \theta) =0, \quad \theta\in[\pi/2,\pi],
\end{equation}

\begin{equation}\label{31b}
L_s^d(0, \theta) =0, \quad \theta\in[0,\pi/2].
\end{equation}

In the basic state, the total intensity $G_s$ is given by:

\begin{equation}\label{32}
G_s=G_s^c+G_s^d=\int_0^{4\pi}[L_s^c(z,\theta)+L_s^d(z,\theta)]d\Omega=L_t\exp\left(-\int_z^1\kappa n_s(z')dz'\right)+\int_0^{\pi}L_s^d(z,\theta)d\Omega,
\end{equation}

where $G_s^c$ and $G_s^d$ represent the collimated and diffuse components of $G_s$, respectively.

Similarly, the radiative heat flux $\boldsymbol{q}_s$ in the basic state can be written as:

\begin{equation}\label{33}
\boldsymbol{q}_s=\boldsymbol{q}_s^c+\boldsymbol{q}_s^d=\int_0^{4\pi}\left(L_s^c(z,\theta)+L_s^d(z,\theta)\right)\boldsymbol{r}d\Omega=L_t\exp\left(\int_z^1-\kappa n_s(z')dz'\right)\hat{\boldsymbol{z}}+\int_0^{4\pi}L_s^d(z,\theta)\boldsymbol{r}d\Omega.
\end{equation}

To solve the problem, a new variable $\tau$ is defined as:

\begin{equation*}
	\tau=\int_z^1 \kappa n_s(z')dz'.
\end{equation*}

These equations provide a set of coupled Fredholm integral equations of the second kind. Specifically, the equation for $G_s$ can be written as:

\begin{equation}\label{34}
G_s(\tau) =e^{-\tau}+ \frac{\omega}{2}\int_0^{\kappa} G_s(\tau')E_1(|\tau-\tau'|)d\tau'+A sgn(\tau-\tau')q_s(\tau')E_2(|\tau-\tau'|),
\end{equation}

where $E_1(x)$ and $E_2(x)$ represent the exponential integrals of order 1 and 2, respectively.

	\begin{equation}\label{35}
	q_s(\tau) = e^{-\tau}+\frac{\omega}{2}\int_0^{\kappa} Aq_s(\tau')E_3(|\tau-\tau'|)d\tau'+sgn(\tau-\tau')G_s(\tau')E_2(|\tau-\tau'|).
\end{equation}
	
The coupled fractional integral equations (FIEs) can be solved using the method of subtraction of singularity, where $E_n(x)$ represents the exponential integral of order $n$ and $sgn(x)$ denotes the signum function.
	
\begin{figure*}[!ht]
	\includegraphics{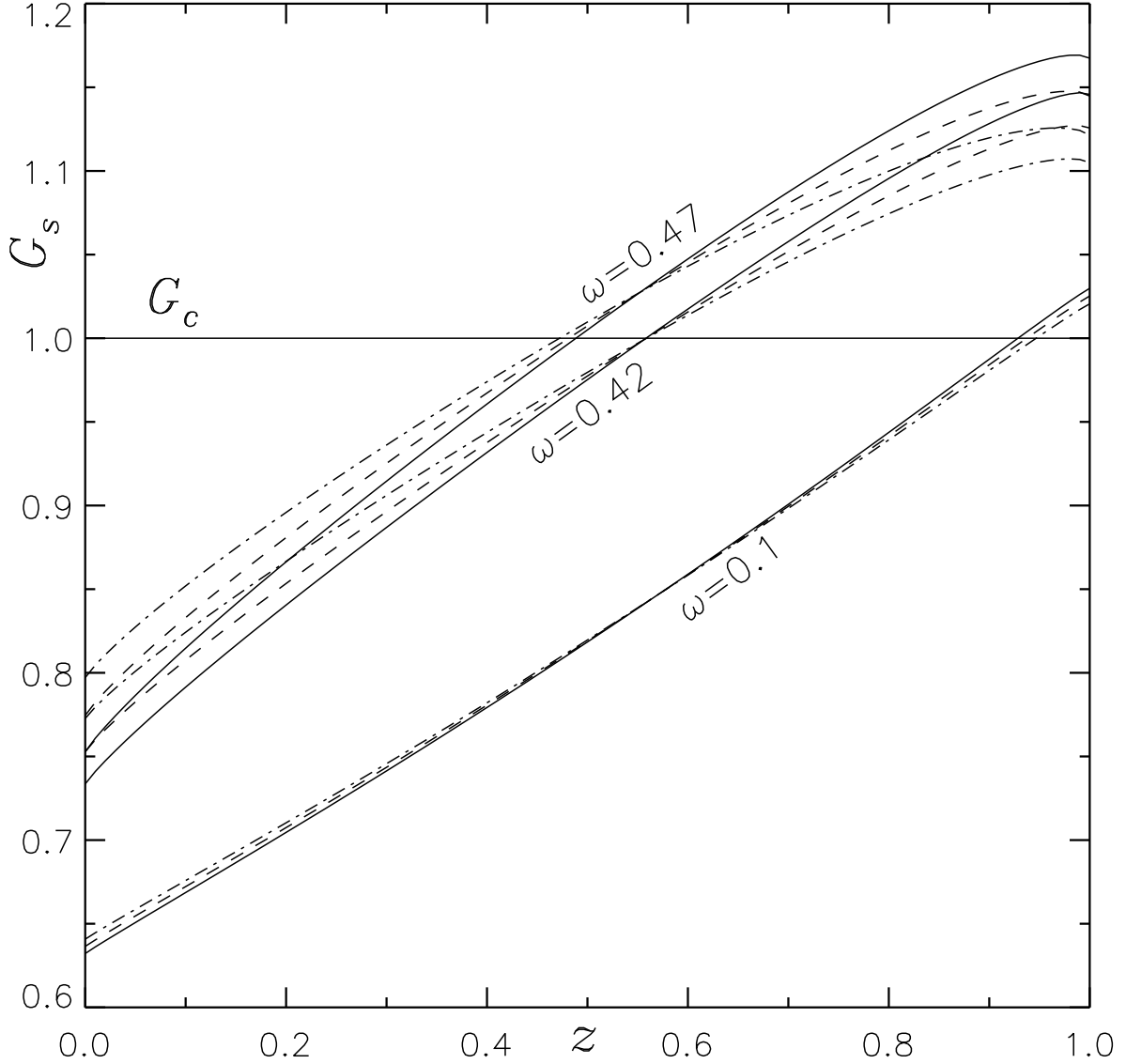}
	\caption{\label{fig2} The change in the overall intensity within a homogeneous suspension by altering the anisotropic scattering coefficient $A$ in the forward direction from 0 to 0.8 is examined for three distinct scenarios: $\omega=0.1,0.42$, and 0.47. The cases are represented as follows: A$=0 (-)$, A$=0.4 (---)$, and A$=0.8 (-\cdot-\cdot-)$. The governing parameter values of the system, namely $S_c=20, k=0.5$, and $L_t=1$, remain constant.}
\end{figure*}
	
The mean swimming direction in the basic state is given by

\begin{equation*}
	<\boldsymbol{p_s}>=-M_s\frac{\boldsymbol{q_s}}{q_s}=M_s\hat{\boldsymbol{z}},
\end{equation*}

where $M_s=M(G_s)$. In the basic state, the cell concentration $n_s(z)$ satisfies the following equation

\begin{equation}\label{36}
\frac{dn_s}{dz}-V_cM_sn_s=0,
\end{equation}

which is accompanied by the equation

\begin{equation}\label{37}
\int_0^1n_s(z)dz=1.
\end{equation}
	
	Equations (\ref{34}) to (\ref{37}) constitute a boundary value problem that can be solved using numerical techniques such as the shooting method.
	
	To demonstrate the equilibrium state, we fix the values of $S_c=20$, $G_c=1$, $\kappa=0.5$, and $\omega$ at specific values of 0.1, 0.42, and 0.47. We then investigate the influence of the forward scattering coefficient $A$ by varying it within the range of 0 to 0.8. The focus is on analyzing the effects of these variations on the total intensity ($G_s$) and the fundamental equilibrium solution.
	
	Figure \ref{fig2} visually illustrates the relationship between the total intensity $G_s$ and the depth of a uniform suspension, while incrementally increasing the forward scattering coefficient $A$ from 0 to 0.8 while keeping $\omega$ constant. As the value of $A$ increases, there is a noticeable increase in the intensity of $G_s$ at lower regions of the uniform suspension. Conversely, at the upper regions, the intensity of $G_s$ decreases.

	\begin{figure*}[!ht]
	\includegraphics{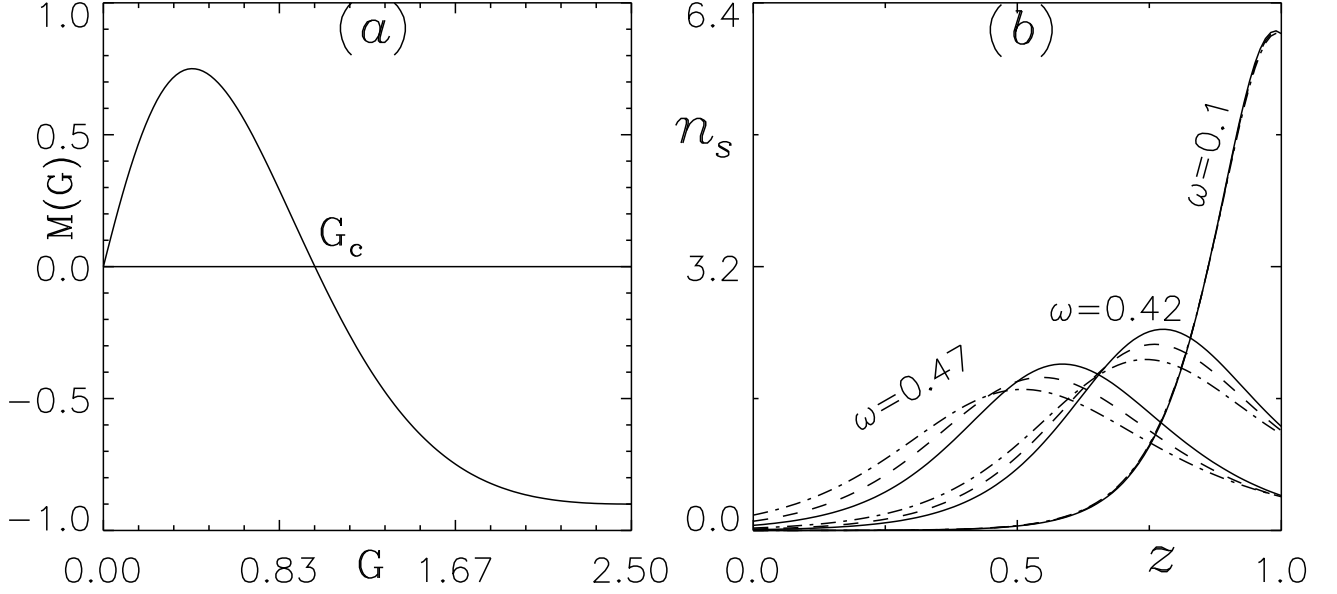}
	\caption{\label{fig3} (a) The response curve of taxis for $G_c=1$, and (b) changes in the concentration profile of the base by varying the forward scattering coefficient $A$ from 0 to 0.8 for three different scenarios: $\omega=0.1,0.42$ and 0.47. The cases are represented as follows: A$=0 (-)$, A$=0.4 (---)$, and A$=0.8 (-\cdot-\cdot-)$. The fixed parameter values governing the system are $S_c=20, V_c=20, k=0.5$, and $L_t=1$.}
\end{figure*}
	
	In Figure \ref{3}(a), the plot illustrates the phototaxis function as a function of light intensity, considering a critical value of $G_c = 1$. On the other hand, Figure \ref{3}(b) displays the influence of the forward scattering coefficient $A$ on the equilibrium position of the sublayer, using the same governing parameters.
	
	In the case where $\omega = 0$, the sublayer forms near the top of the suspension domain in the equilibrium state. However, as $A$ increases from 0 to 0.8, the sublayer's equilibrium position gradually shifts towards the lower region.
	
	Similarly, for $\omega = 0.42$ and 0.47, the sublayer's equilibrium position in the equilibrium state is located approximately at three-quarters and the middle height of the suspension, respectively. Once again, an increase in $A$ from 0 to 0.8 leads to a downward shift of the sublayer's equilibrium position in each scenario.

	\section{Linear stability of the problem}
	To analyse stability, we employ linear perturbation theory, where a small perturbation with an amplitude $\epsilon \ll 1$ is introduced to the equilibrium state. This perturbation can be mathematically represented by the following equation:
	\begin{widetext}
		
		\begin{equation}\label{38}
			[\boldsymbol{U},n,\zeta,L,<p>]=[0,n_s,\zeta_s,L_s,<p_s>]+\epsilon [\boldsymbol{U}_1,n_1,\zeta_1,L_1,<\boldsymbol{p}_1>]+\mathcal{O}(\epsilon^2)  
		\end{equation}
	\end{widetext}
By substituting the perturbed variables into equations (\ref{14}) to (\ref{16}) and linearizing the equations, we gather terms of $o(\epsilon)$ around the equilibrium state, leading to the following expression:
	\begin{equation}\label{39}
		\boldsymbol{\nabla}\cdot \boldsymbol{U}_1=0,
	\end{equation}
	where  $\boldsymbol{U}_1=(U_1,V_1,W_1)$.
\begin{equation}\label{40}
	Sc^{-1}\left(\frac{\partial \boldsymbol{U_1}}{\partial t}\right)+\sqrt{T_a}(z\times U_1)+\boldsymbol{\nabla} P_{e}+Rn_1\hat{\boldsymbol{z}}=\nabla^{2}\boldsymbol{ U_1},
\end{equation}
	\begin{equation}\label{41}
		\frac{\partial{n_1}}{\partial{t}}+V_c\boldsymbol{\nabla}\cdot(<\boldsymbol{p_s}>n_1+<\boldsymbol{p_1}>n_s)+W_1\frac{dn_s}{dz}=\boldsymbol{\nabla}^2n_1.
	\end{equation}
	If $G=G_s+\epsilon G_1+\mathcal{O}(\epsilon^2)=(G_s^c+\epsilon G_1^c)+(G_s^d+\epsilon G_1^d)+\mathcal{O}(\epsilon^2)$, then the steady collimated total intensity is perturbed as $L_t\exp\left(-\kappa\int_z^1(n_s(z')+\epsilon n_1+\mathcal{O}(\epsilon^2))dz'\right)$  and after simplification, we get
	\begin{equation}\label{42}
		G_1^c=L_t\exp\left(-\int_z^1 \kappa n_s(z')dz'\right)\left(\int_1^z\kappa n_1 dz'\right)
	\end{equation}
	Similarly, $G_1^d$ is given by 
	\begin{equation}\label{43}
		G_1^d=\int_0^{4\pi}L_1^d(\boldsymbol{ x},\boldsymbol{ r})d\Omega.
	\end{equation}
	Similarly, for the radiative heat flux $q=q_s+\epsilon q_1++\mathcal{O}(\epsilon^2)=(q_s^c+\epsilon q_1^c)+(q_s^d+\epsilon q_1^d)+\mathcal{O}(\epsilon^2)$, we find 
	\begin{equation}\label{44}
		\boldsymbol{q}_1^c=L_t\exp\left(-\int_z^1 \kappa n_s(z')dz'\right)\left(\int_1^z\kappa n_1 dz'\right)\hat{z}
	\end{equation}
	and
	\begin{equation}\label{45}
		q_1^d=\int_0^{4\pi}L_1^d(\boldsymbol{ x},\boldsymbol{r})\boldsymbol{ r}d\Omega.
	\end{equation}
	After perturbing the expression for swimming orientation and collecting $O(\epsilon)$ terms gives the perturbed swimming orientation
	\begin{equation}\label{46}
		<\boldsymbol{p_1}>=G_1\frac{dM_s}{dG}\hat{\boldsymbol{z}}-M_s\frac{\boldsymbol{q_1}^H}{\boldsymbol{q_s}},
	\end{equation}
here, $\boldsymbol{q}_1^H=[\boldsymbol{q}_1^x,\boldsymbol{q}_1^y]$ represents the horizontal component of the perturbed radiative flux $\boldsymbol{q}_1$.
By substituting the value of $<\boldsymbol{p_1}>$ from Eq.$(\ref{46})$ into Eq.$(\ref{41})$ and simplifying, we obtain:
	\begin{equation}\label{47}
		\frac{\partial{n_1}}{\partial{t}}+V_c\frac{\partial}{\partial z}\left(M_sn_1+n_sG_1\frac{dM_s}{dG}\right)-V_cn_s\frac{M_s}{q_s}\left(\frac{\partial q_1^x}{\partial x}+\frac{\partial q_1^y}{\partial y}\right)+W_1\frac{dn_s}{dz}=\nabla^2n_1.
	\end{equation}
	By eliminating $P_e$ and the horizontal component of $u_1$, equations (\ref{40}) and (\ref{41}) can be simplified to three equations for the perturbed variables: the vertical component of the velocity $w_1$, the vertical component of the vorticity $\zeta_1$ (which is defined as $\zeta\cdot\hat{\boldsymbol{z}}$), and the concentration $n_1$. These variables can be further decomposed into normal modes as:
	
	\begin{equation}\label{48}
	[W_1,\zeta_1,n_1]=[\tilde{W}(z),\tilde{Z}(z),\tilde{N}(z)]\exp{(\sigma t+i(lx+my))}. 
	\end{equation}
The equation governing the perturbed intensity $L_1$ can be written as:
	\begin{equation}\label{49}
		\xi\frac{\partial L_1}{\partial x}+\eta\frac{\partial L_1}{\partial y}+\nu\frac{\partial L_1}{\partial z}+\kappa( n_sL_1+n_1L_s)=\frac{\omega\kappa}{4\pi}(n_sM_1+G_sn_1+A\nu(n_sq_s\cdot\hat{z}-q_sn_1))-\kappa L_sn_1,
	\end{equation}
	with boundary conditions
	\begin{subequations}
		\begin{equation}\label{50a}
			L_1(x, y, z = 1, \xi, \eta, \nu) =0,\qquad \theta\in[\pi/2,\pi], ~~\phi\in[0, 2\pi], 
		\end{equation}
		\begin{equation}\label{50b}
			L_1(x, y, z = 0,\xi, \eta, \nu) =0,\qquad \theta\in[0,\pi/2], ~~\phi\in[0, 2\pi]. 
		\end{equation}
	\end{subequations}
	The Eq. $(\ref{48})$ implies that $L_1^d$ can be represented by the following expression:
	\begin{equation*}
		L_1^d=\Psi^d(z,\xi,\eta,\nu)\exp{(\sigma t+i(lx+my))}. 
	\end{equation*}
	From Eqs.~(\ref{42}) and (\ref{43}), we get
	\begin{equation}\label{51}
		G_1^c=\left[L_t\exp\left(-\int_z^1 \kappa n_s(z')dz'\right)\left(\int_1^z\kappa n_1 dz'\right)\right]\exp{(\sigma t+i(lx+my))}=\mathbb{G}^c(z)\exp{(\sigma t+i(lx+my))},
	\end{equation}
	and 
	\begin{equation}\label{52}
		G_1^d=\mathbb{G}^d(z)\exp{(\sigma t+i(lx+my))}= \left(\int_0^{4\pi}\Psi^d(z,\xi,\eta,\nu)d\Omega\right)\exp{(\sigma t+i(lx+my))},
	\end{equation}

	where $\mathbb{G}(z)=\mathbb{G}^c(z)+\mathbb{G}^d(z)$ is the perturbed total intensity. Similarly from Eqs.~(\ref{44}) and (\ref{45}), we have
	\begin{equation*}
		\boldsymbol{q}_1=[q_1^x,q_1^y,q_1^z]=[P(z),Q(z),S(z)]\exp{[\sigma t+i(lx+my)]},
	\end{equation*}
	where
	\begin{equation*}
		[P(z), Q(z)]=\int_0^{4\pi}[\xi,\eta]\Psi^d(z,\xi,\eta,\nu) d\Omega,
	\end{equation*}
and
    \begin{equation*}
    	S(z)=\int_0^{4\pi}[\Psi^c(z,\xi,\eta,\nu)+\Psi^d(z,\xi,\eta,\nu)]\nu d\Omega.
    \end{equation*}
 
	Note that the $P(z)$ and $Q(z)$ appear due to scattering. On the other hand, $S(z)$ appears due to anisotropic scattering which becomes zero in the case of isotropic scattering.
	
	Now $\Psi^d$ satisfies
	\begin{equation}\label{53}
		\frac{d\Psi^d}{dz}+\frac{(i(l\xi+m\eta)+\kappa n_s)}{\nu}\Psi^d=\frac{\omega\kappa}{4\pi\nu}(n_s\mathcal{\mathbb{G}}+G_s\tilde{N}(z)+A\nu(n_s S-q_s\tilde{N}(z)))-\frac{\kappa}{\nu}I_s\tilde{N}(z),
	\end{equation}
	subject to the boundary conditions
	\begin{subequations}
		\begin{equation}\label{54a}
			\Psi^d( 1, \xi, \eta, \nu) =0,\qquad \theta\in[\pi/2,\pi], ~~\phi\in[0, 2\pi] , 
		\end{equation}
		\begin{equation}\label{54b}
			\Psi^d( 0,\xi, \eta, \nu) =0,\qquad \theta\in[0,\pi/2], ~~\phi\in[0, 2\pi]. 
		\end{equation}
	\end{subequations}

	The stability equations reformed as
	\begin{equation}\label{55}
		\left(\sigma S_c^{-1}+k^2-\frac{d^2}{dz^2}\right)\left( \frac{d^2}{dz^2}-k^2\right)\tilde{W}=Rk^2\tilde{N}(z),
	\end{equation}
	\begin{equation}\label{56}
		\left(\sigma S_c^{-1}+k^2-\frac{d^2}{dz^2}\right)\tilde{Z}(z)=\sqrt{T_a}\frac{d\tilde{W}}{dz}
	\end{equation}
	\begin{equation}\label{57}
		\left(\sigma+k^2-\frac{d^2}{dz^2}\right)\tilde{N}(z)+V_c\frac{d}{dz}\left(M_s\tilde{N}(z)+n_s\mathbb{G}\frac{dM_s}{dG}\right)-i\frac{V_cn_sM_s}{q_s}(lP+mQ)=-\frac{dn_s}{dz}\tilde{W}(z),
	\end{equation}
	
	with boundary conditions
	\begin{equation}\label{58}
	\tilde{W}(z)=\frac{d\tilde{W}(z)}{dz}=\tilde{Z}(z)=\frac{d\tilde{N}(z)}{dz}-V_cM_s\tilde{N}(z)-n_sV_c\mathbb{G}\frac{dM_s}{dG}=0\quad at\quad z=0,
	\end{equation}
	\begin{equation}\label{59}
	\tilde{W}(z)=\frac{d\tilde{W}(z)}{dz}=\frac{d\tilde{Z}(z)}{dz}=\frac{d\tilde{N}(z)}{dz}-V_cM_s\tilde{N}(z)-n_sV_c\mathbb{G}\frac{dM_s}{dG}=0\quad at\quad z=1.
	\end{equation}
	
n the given context, the non-dimensional wavenumber is represented by $k$, which is determined by taking the square root of the sum of the squares of $l$ and $m$. Equations (\ref{58})-(\ref{59}) form an eigenvalue problem, where $\sigma$ is described as a function of several dimensionless parameters such as $V_c$, $\kappa$, $\omega$, $A$, $l$, $m$, and $R$. Equation (\ref{57}) can be equivalently expressed as follows:
	\begin{equation}\label{60}
       D^2\tilde{N}(z)-\Lambda_3(z)D\tilde{N}(z)-(\sigma+k^2+\Lambda_2(z))\tilde{N}(z)-\Lambda_1(z)\int_1^z\tilde{N}(z) dz-\Lambda_0(z)=Dn_s\tilde{W}, 
	\end{equation}
	where
	\begin{subequations}
		\begin{equation}\label{61a}
			\Lambda_0(z)=V_cD\left(n_s\mathbb{G}^d\frac{dM_s}{dG}\right)-\iota\frac{V_cn_sM_s}{q_s}(lP+mQ),
		\end{equation}
		\begin{equation}\label{61b}
			\Lambda(z)=\kappa V_cD\left(n_sG_s^c\frac{dM_s}{dG}\right)
		\end{equation}
		\begin{equation}\label{61c}
			\Lambda_2(z)=2\kappa V_c n_s G_s^c\frac{dM_s}{dG}+V_c\frac{dM_s}{dM}DG_s^d,
		\end{equation}
		\begin{equation}\label{61d}
			\Lambda_3(z)=V_cM_s.
		\end{equation}
	\end{subequations}
	Introducing a novel variable denoted as
	\begin{equation}\label{62}
		\tilde{\Theta}(z)=\int_1^z\tilde{N}(z')dz',
	\end{equation}
 the linear stability equations can be reformulated as
	\begin{equation}\label{63}
		D^4\tilde{W}-\left(\sigma S_c^{-1}+k^2\right)D^2\tilde{W}-\left(\sigma S_c^{-1}+k^2\right)\tilde{W}=Rk^2D\tilde{\Theta},
	\end{equation}
	\begin{equation}\label{64}
		D^2\tilde{Z}(z)-\left(\sigma S_c^{-1}+k^2-D^2\right)\tilde{Z}(z)=\sqrt{T_a}D\tilde{W}
	\end{equation}
	\begin{equation}\label{65}
	D^3\tilde{\Theta}-\Lambda_3(z)D^2\tilde{\Theta}-(\sigma+k^2+\Lambda_2(z))D\tilde{\Theta}-\Lambda_1(z)\tilde{\Theta}-\Lambda_0(z)=Dn_s\tilde{W}. 
	\end{equation}
The boundary conditions remain unchanged except for the flux boundary condition, given by the equation

\begin{equation}\label{66}
D^2\tilde{N}(z)-\Lambda_2D\tilde{N}(z)-n_sV_c\mathbb{G}\frac{dM_s}{dG}=0\quad at\quad z=0,1.
\end{equation}
Furthermore, an additional boundary condition is introduced
	\begin{equation}\label{67}
		\tilde{\Theta}(z)=0,\quad at\quad z=1.
	\end{equation}

	\section{SOLUTION PROCEDURE}

To obtain the neutral (marginal) stability curves or the growth rate $Re(\sigma)$ as a function of R in the (k, R)-plane for a fixed parameter set, a fourth-order accurate finite-difference scheme based on Newton-Raphson-Kantorovich (NRK) iterations is utilized to solve equations \ref{63} to \ref{65}. By analysing the resulting graph, points where $Re(\sigma) = 0$ are identified, representing the marginal (neutral) stability curve. On this curve, if $Im(\sigma) = 0$, it indicates a stationary (non-oscillatory) bioconvective solution, while non-zero $Im(\sigma)$ indicates oscillatory solutions. If the most unstable mode remains on the oscillatory branch of the neutral curve, it signifies overstability.

When oscillatory solutions are present, there is a common point $(k_b)$ where the stationary and oscillatory branches intersect, and the oscillatory branch includes points with $k \leq k_b$. Among the points on the neutral curve $R^{(n)}(k)$ (where n = 1, 2, 3, ...), a particular most unstable mode $(k_c, R_c)$ is selected. The wavelength of the initial disturbance is then determined as $\lambda_c = 2\pi/k_c$. The bioconvective solution is referred to as mode n if it can be organized into n convection cells such that one cell overlaps another vertically.

	\section{NUMERICAL RESULTS}
To identify the most unstable mode starting from an initial equilibrium solution, we employ a specific set of fixed parameters. These parameters include $S_c = 20$, $G_c = 1$, $V_c = 10$, $15$, $20$, with $\omega$ ranging from $0$ to $1$, $\kappa$ values of $0.5$ and $1$, and $A$ values of $0$, $0.4$, and $0.8$. The Taylor number, which represents the rotation rate, is varied from lower non-zero values to higher values. To investigate the impact of rotation on a forward scattering algal suspension, we focus on specific values of the scattering albedo $\omega$. For the case of $\kappa=0.5$, we examine the scenarios where $\omega$ takes on the values $0.1$, $0.42$, and $0.47$. Similarly, for the case of $\kappa=1$, we consider $\omega$ values of $0.1$, $0.58$, and $0.60$. When $\omega=0.1$, the equilibrium state's sublayer is positioned near the top of the domain. As the value of $\omega$ increases, the sublayer's equilibrium position shifts from the top to approximately three-quarters height, and then further to the middle height of the domain.
	
	\subsection{WHEN SCATTERING IS WEAK  }
The primary focus of this study is to investigate how rotation influences the initiation of bioconvection. To explore this, the researchers compare the relative effectiveness of self-shading and scattering, specifically by utilizing a lower scattering albedo value denoted as $\omega$. Additionally, the strength of self-shading is manipulated by varying the extinction coefficient $\kappa$, with a higher value of $\kappa=1$ indicating strong self-shading and a lower value of $\kappa=0.5$ representing weak self-shading. Throughout the study, a consistent critical intensity of $G_c=1$ is used. The researchers divide the entire study into three distinct cases based on the position of the equilibrium state, namely at the top, three-quarter, and mid-height of the domain. Suitable values of $\omega$ are chosen accordingly for each case.

\subsubsection{$V_c=20$}
	(a) $\kappa=0.5$
	
	This section delves into the examination of bioconvective instability and its correlation with the Taylor number. The study investigates this relationship by considering three distinct values of the forward scattering coefficient $A$ and three different values of $\omega$, while keeping the parameters $V_c$ and $\kappa$ fixed at 20 and 0.5, respectively. Similar to previous sections, the study is divided into three cases, each based on the position of the equilibrium state within the domain (at the top, three-quarters, and mid-height). The appropriate value of $\omega$ is chosen accordingly for each case.
	
\begin{figure*}[!htbp]
		\includegraphics{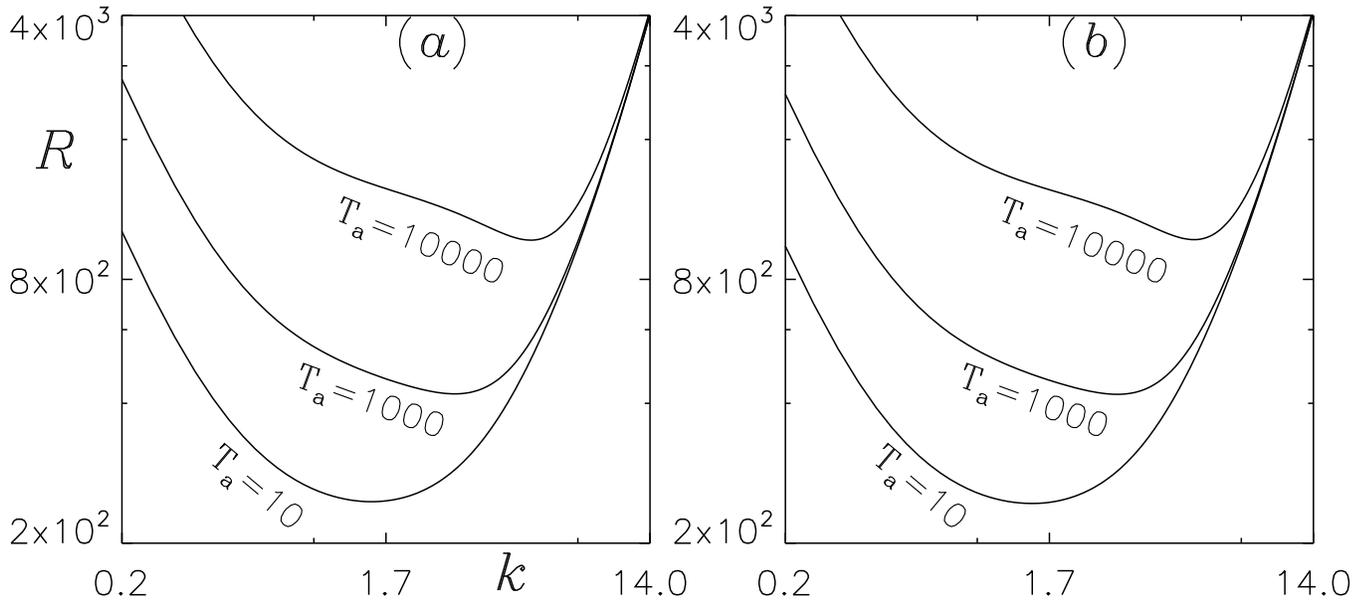}
		\caption{\label{fig4} The marginal stability curves for (a) $A=0$, (b) $A=0.4$. Here, the other governing parameter values $S_c=20,V_c=20,k=0.5$, and $\omega=0.1$ are kept fixed.}
\end{figure*}
	
When $\omega=0.1$, the sublayer is located approximately at the top of the domain when $A=0$. However, as the value of $A$ increases, the sublayer at the equilibrium state shifts towards the three-quarter height of the domain. In this scenario, we vary the Taylor number ($T_a$) from a lower value (e.g., $T_a=10$) to a higher value (e.g., $T_a=10000$). At $T_a=10$, a stationary mode of the bioconvective solution is observed. As $T_a$ is increased up to 10000, the behaviour of the solutions remains the same, but the critical Rayleigh number increases (refer to Fig.\ref{fig4}(a)). For $A=0.4$ and $A=0.8$, the mode of the bioconvective solution remains unchanged, and as $T_a$ is increased, the critical Rayleigh number also increases (refer to Fig.\ref{fig4}).
	
\begin{figure*}[!bt]
	\includegraphics{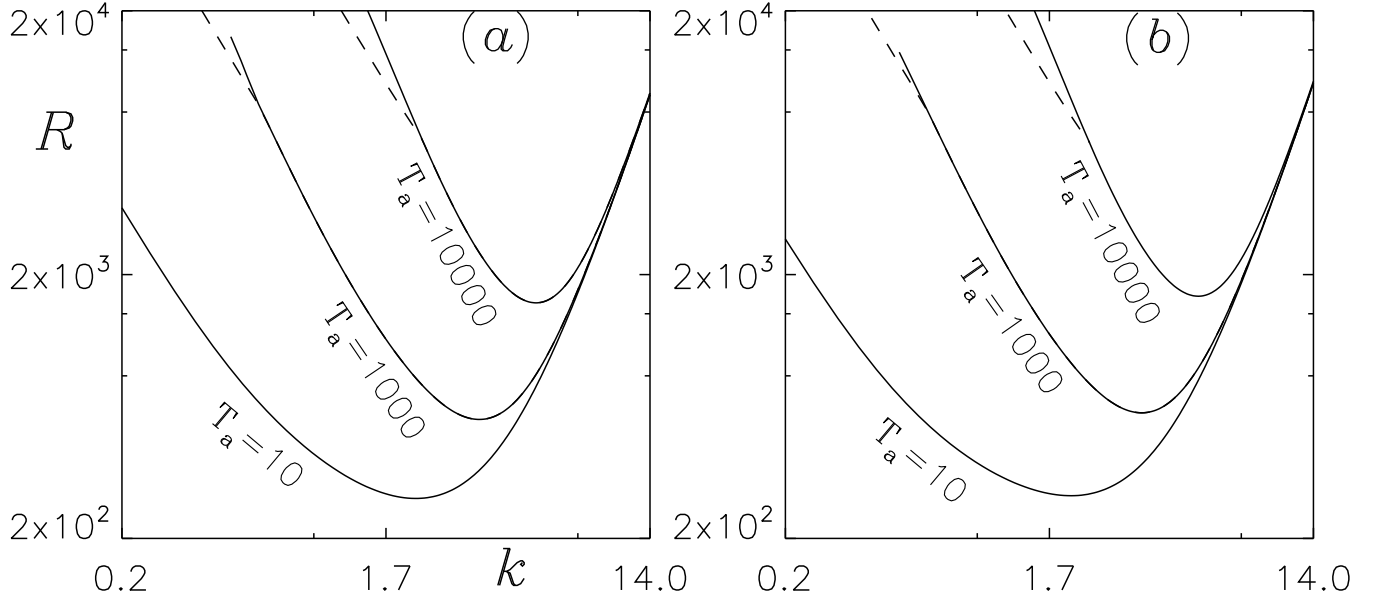}
	\caption{\label{fig5}The marginal stability curves for (a) $A=0$, (b) $A=0.4$. Here, the other governing parameter values $S_c=20,V_c=20,k=0.5$, and $\omega=0.42$ are kept fixed.}
\end{figure*}	
	
When $\omega=0.42$, the sublayer in the equilibrium state is located at the three-quarter height of the suspension for $A=0$. At $T_a=10$, a stationary mode of the bioconvective solution is observed. As $T_a$ is increased to 1000, an oscillatory branch splits from the stationary branch of the marginal stability curve. However, the most unstable solution still occurs on the stationary branch of the neutral curve, resulting in a stationary solution and an increased critical Rayleigh number. Similarly, at $T_a=10000$, the same nature of the marginal stability curve and a bioconvective solution is observed, with an increased critical Rayleigh number. As $A$ increases, the sublayer at the equilibrium state shifts towards the mid-height of the domain. For $A=0.4$ and $A=0.8$, the mode of the bioconvective solution remains unchanged, and as $T_a$ is increased, the critical Rayleigh number also increases (refer to Fig.~\ref{fig5}).

\begin{figure*}[!htbp]
	\includegraphics{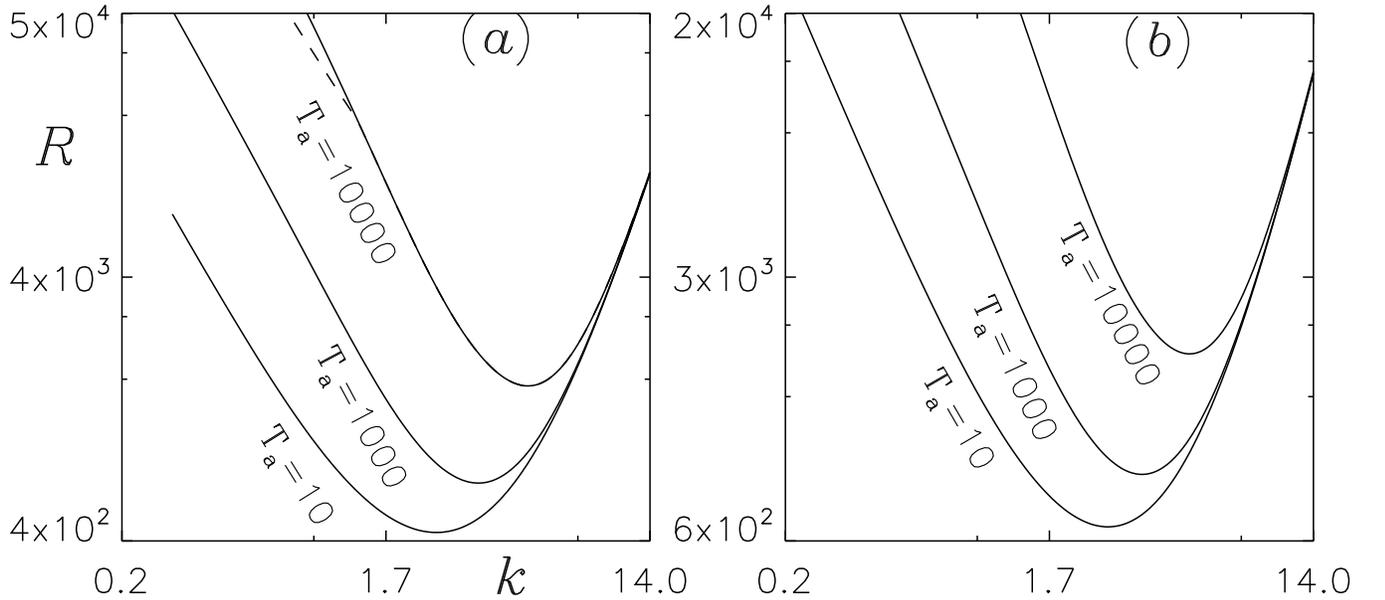}
	\caption{\label{fig6}The marginal stability curves for (a) $A=0$, (b) $A=0.4$. Here, the other governing parameter values $S_c=20,V_c=20,k=0.5$, and $\omega=0.47$ are kept fixed.}
\end{figure*}
When $\omega=0.47$, the sublayer in the equilibrium state is located approximately at the middle height of the domain for $A=0$. At $T_a=10$, a stationary mode of the bioconvective solution is observed. As $T_a$ is increased to 1000, the mode of the bioconvective instability remains the same, but the critical Rayleigh number increases. As $T_a$ reaches 10000, an oscillatory branch bifurcates from the stationary branch of the marginal stability curve at $k\approx1.6$ and exists for all $k<1.6$. However, the most unstable bioconvective solution still occurs on the stationary branch, resulting in a stationary solution. As $A$ increases, the sublayer at the equilibrium state shifts towards the bottom of the domain. For $A=0.4$, the stationary nature of the solution is observed at $T_a=10$, and this nature of the solution remains unchanged for all higher values of the Taylor number. However, the critical Rayleigh number increases as $T_a$ is increased. A similar nature of the bioconvective solution is also observed for $A=0.8$ (see Fig.~\ref{fig6}).

The numerical results of this section are summarized in Table~\ref{tab1}.

\begin{table}[!htbp]
	\caption{\label{tab1}The table shows the numerical results of bioconvective solutions for a constant $V_c=$20 and $\kappa=$0.5. The values are presented for different increments in $T_a$, while keeping all other parameters unchanged. The presence of a star symbol indicates that the $R^{(1)}(k)$ branch of the neutral curve exhibits oscillatory behaviour.}
	\begin{ruledtabular}
		\begin{tabular}{cccccccccccc}
			\multirow{2}{*}{$\omega$}&\multirow{2}{*}{$T_a$}&\multicolumn{3}{c}{$A=0$} & \multicolumn{3}{c}{$A=0.4$} & \multicolumn{3}{c}{$A=0.8$}\\\cline{3-5}\cline{6-8}\cline{9-11} & & $\lambda_c$ & $R_c$ & $Im(\sigma)$ &$ \lambda_c$ & $R_c$ & $Im(\sigma)$ &$\lambda_c$ & $R_c$ & $Im(\sigma)$\\
			\hline
			
			&0     & 4.23 & 194.18 & 0 & 4.23 & 193.15 & 0 & 4.39 & 192.07 & 0 \\ 		
			0.1&1000  & 2.17 & 379.88 & 0 & 2.17 & 379.37 & 0 & 2.25 & 379.09 & 0 \\ 		
			&10000 & 1.22 & 988.77 & 0 & 1.17 & 990.33 & 0 & 1.18 & 991.91 & 0 \\ 
			\hline
			
			&0     & 2.98 & 217.27  & 0 & 3.15 & 222.85  & 0 & 3.44 & 231.89 & 0 \\ 		
			0.42&1000  & 1.17$^{\star}$ & 452.01  & 0 & 1.79$^{\star}$ & 480.72  & 0 & 1.80 & 522.28 & 0 \\ 		
			&10000 & 1.12$^{\star}$ & 1332.85 & 0 & 1.13$^{\star}$ & 1417.52 & 0 & 1.14 & 1533.90 & 0 \\ 	
			\hline
			
			&0     & 2.53 & 432.26  & 0 & 2.35 & 616.06  & 0 & 2.17 & 919.01 & 0 \\ 		
			0.47&1000  & 1.78 & 679.23  & 0 & 1.79 & 879.18  & 0 & 1.76 & 1207.16 & 0 \\ 		
			&10000 & 1.20$^{\star}$ & 1652.82 & 0 & 1.21 & 1989.33 & 0 & 1.23 & 2524.01 & 0 \\  		
		\end{tabular}
	\end{ruledtabular}
	
\end{table}

	(b) $\kappa=1$
In this section, we examine the influence of the Taylor number on bioconvective instability for three different values of the forward scattering coefficient $A$ at three different values of $\omega$. The investigation is conducted with a set of fixed parameters, namely $V_c=20$ and $\kappa=1$.

\begin{figure*}[!htbp]
	\includegraphics{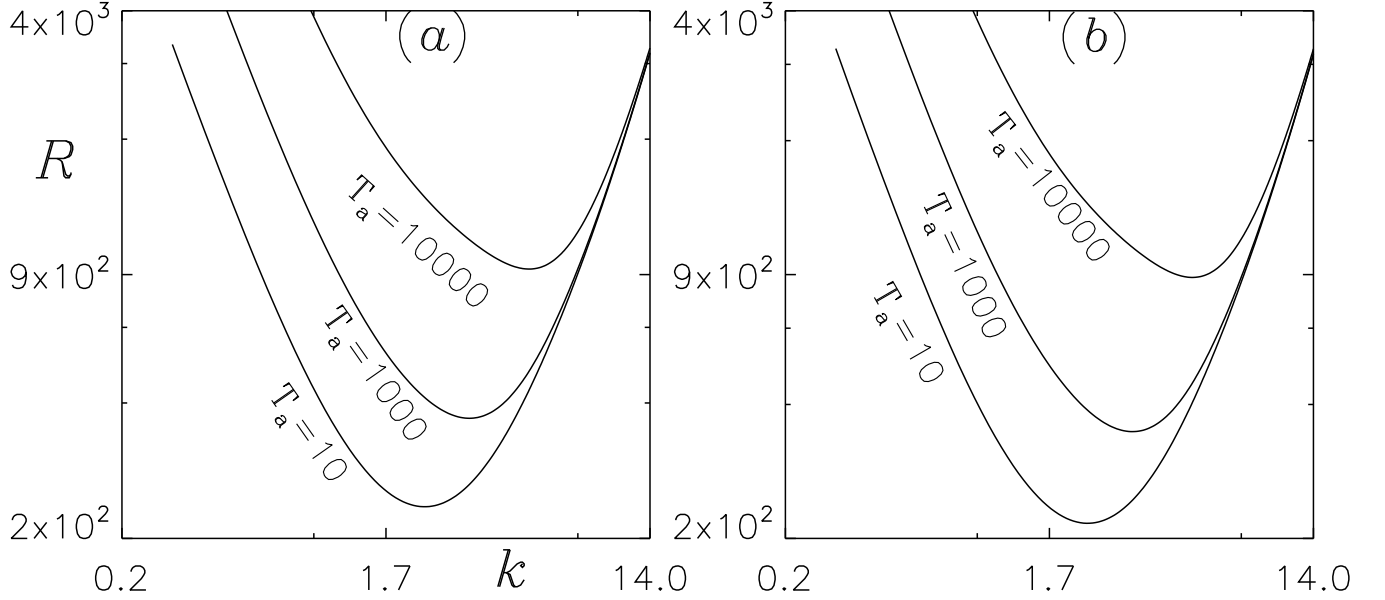}
	\caption{\label{fig7} The marginal stability curves for (a) $A=0$, (b) $A=0.4$. Here, the other governing parameter values $S_c=20,V_c=20,k=1$, and $\omega=0.1$ are kept fixed.}
\end{figure*} 

For $\omega=0.1$, the sublayer occurs approximately at the top of the domain when $A=0$. As the value of $A$ increases, the sublayer at the equilibrium state shifts towards the three-quarter height of the domain. In this case, we vary the Taylor number from a lower value of $T_a=10$ to a higher value of $T_a=10000$. At $T_a=10$, a stationary mode of the bio-convective solution is observed, and as $T_a$ is increased up to 10000, the same behaviour of solutions is observed with an increase in the critical Rayleigh number (refer to Fig.~\ref{fig7}(a)). For $A=0.4$ and $A=0.8$, the mode of the bio-convective solution remains unchanged, and as $T_a$ is increased, the critical Rayleigh number also increases (see Fig.~\ref{fig4}).

\begin{figure*}[!htbp]
	\includegraphics{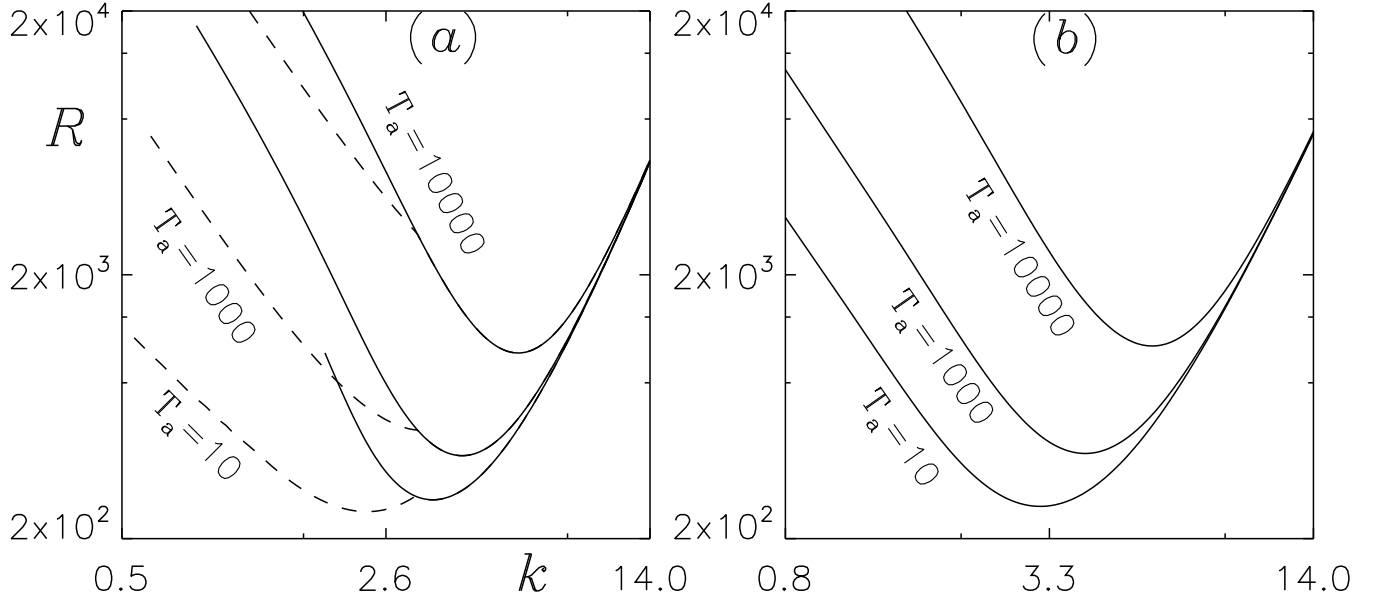}
	\caption{\label{fig8} The marginal stability curves for (a) $A=0$, (b) $A=0.8$. Here, the other governing parameter values $S_c=20,V_c=20,k=1$, and $\omega=0.58$ are kept fixed.}
\end{figure*}

For $\omega=0.58$, the sublayer in the equilibrium state occurs at the three-quarter height of the suspension when $A=0$. At $T_a=10$, an oscillatory branch splits from the stationary branch of the marginal stability curve at $k\approx3$ and exists for all $k<3$. In this case, the most unstable solution occurs on the oscillatory branch of the marginal stability curve, resulting in an overstable solution. As $T_a$ is increased to 1000, another oscillatory branch splits from the stationary branch, but the most unstable solution occurs on the stationary branch of the neutral curve, leading to a stationary solution and an increase in the critical Rayleigh number. For $T_a=10000$, the same nature of the marginal stability curve and a bioconvective solution is observed, and the critical Rayleigh number continues to increase (refer to Fig.~\ref{fig8}(a)). As the value of $A$ increases, the sublayer at the equilibrium state shifts towards the mid-height of the domain. For $A=0.4$, the nature of the bioconvective solution and the critical Rayleigh number remain unchanged compared to the case when $A=0$. For $A=0.8$ and $T_a=0.4$, a stationary nature of the bioconvective instability is observed for all values of the Taylor number, and the critical Rayleigh number increases as $T_a$ is increased, similar to the previous case (see Fig.~\ref{fig8}).

\begin{figure*}[!htbp]
	\includegraphics{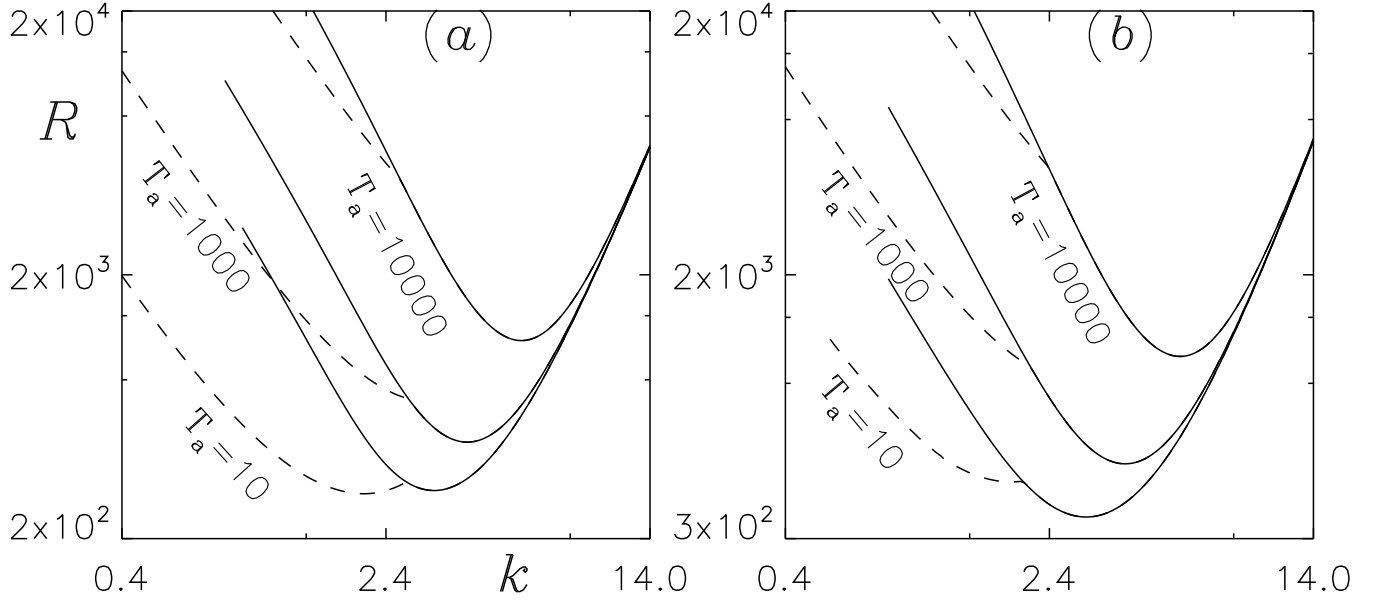}
	\caption{\label{fig9} The marginal stability curves for (a) $A=0$, (b) $A=0.4$. Here, the other governing parameter values $S_c=20,V_c=20,k=1$, and $\omega=0.6$ are kept fixed.}
\end{figure*}

For $\omega=0.6$, the sublayer in the equilibrium state occurs at the three-quarter height of the suspension when $A=0$. At $T_a=10$, an oscillatory branch splits from the stationary branch of the marginal stability curve at $k\approx3$ and exists for all $k<3$. In this case, the most unstable solution occurs on the oscillatory branch of the marginal stability curve, resulting in an overstable solution. As $T_a$ is increased to 1000, another oscillatory branch splits from the stationary branch, but the most unstable solution occurs on the stationary branch of the neutral curve, leading to a stationary solution and an increase in the critical Rayleigh number. For $T_a=10000$, the same nature of the marginal stability curve and a bioconvective solution is observed, and the critical Rayleigh number continues to increase (refer to Fig.~\ref{fig9}(a)). As the value of $A$ increases, the sublayer at the equilibrium state shifts towards the mid-height of the domain. For $A=0.4$ and $T_a=10$, an oscillatory branch splits from the stationary branch of the marginal stability curve at $k\approx2.2$ and exists for all $k<2.2$. Here, the most unstable solution occurs on the stationary branch of the marginal stability curve, resulting in a stationary solution. As $T_a$ is increased to 1000 and, 10000, another oscillatory branch splits from the stationary branch, but the most unstable solution occurs on the stationary branch of the neutral curve, leading to a stationary solution and an increase in the critical Rayleigh number (refer to Fig.~\ref{fig9}(b)). For $T_a=0.8$, a stationary nature of the bioconvective instability is observed for all values of the Taylor number, and the critical Rayleigh number increases as $T_a$ is increased, similar to the previous cases (refer to Table~\ref{tab2}).

The numerical results of this section are presented in Table~\ref{tab2}.

\begin{table}[!htbp]
	\caption{\label{tab2}The table shows the numerical results of bioconvective solutions for a constant $V_c=$20 and $\kappa=$0.5. The values are presented for different increments in $T_a$, while keeping all other parameters unchanged. The presence of a star symbol indicates that the $R^{(1)}(k)$ branch of the neutral curve exhibits oscillatory behaviour.}
	\begin{ruledtabular}
		\begin{tabular}{cccccccccccc}
			\multirow{2}{*}{$\omega$}&\multirow{2}{*}{$T_a$}&\multicolumn{3}{c}{$A=0$} & \multicolumn{3}{c}{$A=0.4$} & \multicolumn{3}{c}{$A=0.8$}\\\cline{3-5}\cline{6-8}\cline{9-11} & & $\lambda_c$ & $R_c$ & $Im(\sigma)$ &$ \lambda_c$ & $R_c$ & $Im(\sigma)$ &$\lambda_c$ & $R_c$ & $Im(\sigma)$\\
			\hline
			
			&0     & 2.78 & 239.43 & 0 & 2.78 & 239.06  & 0    & 2.76 & 238.56 & 0 \\ 		
		 0.1&1000  & 1.92 & 395.58 & 0 & 1.92 & 395.67  & 0    & 1.92 & 395.64 & 0 \\ 		
			&10000 & 1.18 & 922.98 & 0 & 1.19 & 924.11  & 0    & 1.19 & 924.11 & 0 \\ 
			\hline
			
			&0     & 2.69$^{\ddag}$ & 312.75$^{\ddag}$  & 8.79 & 2.82$^{\ddag}$ & 304.94$^{\ddag}$  & 7.65 & 1.98 & 326.80 & 0 \\ 		
		0.58&1000  & 1.48$^{\star}$ & 498.04  & 0    & 1.50$^{\star}$ & 496.53  & 0    & 1.55 & 507.05 & 0 \\ 		
			&10000 & 1.03$^{\star}$ & 1168.45 & 0    & 1.05$^{\star}$ & 1193.34 & 0    & 1.08 & 1238.14 & 0 \\ 	
			\hline
			
			&0     & 3.09$^{\ddag}$ & 295.61$^{\ddag}$  & 7.16 & 2.09$^{\star}$ & 310.62  & 0    & 2.05 & 475.61 & 0 \\ 		
		 0.6&1000  & 1.54$^{\star}$ & 464.30  & 0    & 1.61$^{\star}$ & 481.28  & 0    & 1.66 & 644.02 & 0 \\ 		
			&10000 & 1.07$^{\star}$ & 1126.63 & 0    & 1.11 & 1164.71 & 0    & 1.15 & 1376.95 & 0 \\  		
		\end{tabular}
	\end{ruledtabular}
	
\end{table}	
	
%
	
	\section{Conclusion}
The proposed phototaxis model investigates the impact of rotation on the initiation of light-induced bioconvection in a forward-scattering algal suspension illuminated by collimated flux. The model incorporates a linear anisotropic scattering coefficient.

When the forward scattering coefficient $A$ is increased in a uniform suspension, several effects are observed. Firstly, the total intensity of light undergoes changes in the equilibrium state. Specifically, it decreases in the upper half of the suspension and increases in the lower half. This indicates a redistribution of light intensity within the suspension. Furthermore, the critical value of total intensity, which is the threshold for the onset of bioconvection, exhibits variations. In the lower half of the suspension, the critical intensity decreases as $A$ is increased, implying that bioconvection is more likely to occur at lower light intensities. Conversely, in the upper half of the suspension, the critical intensity increases with increasing $A$, indicating a higher threshold for bioconvection initiation. Additionally, the position of the sublayer, where cells accumulate, is influenced by the value of $A$. As $A$ is increased, the sublayer at the equilibrium state shifts towards the bottom of the domain. This suggests that the presence of forward scattering promotes the accumulation of cells in the lower region of the suspension.

The linear stability analysis predicts both stationary and oscillatory disturbances. The nature of the oscillatory bioconvective solutions appears to depend on the location of the sublayer at the equilibrium state. When the sublayer forms near the top of the domain, the bioconvective instability remains stationary for all values of the Taylor number $T_A$ and forward scattering coefficient $A$. However, when the sublayer forms at the three-quarter height and mid-height of the domain, an overstable nature of the solution is observed for lower values of the Taylor number and forward scattering coefficient.

Furthermore, it is observed that the critical wavelength decreases and the corresponding Rayleigh number increases as the Taylor number increases. This suggests that the suspension becomes more stable as the rotation rate increases.

The proposed model appears to be more realistic than the model presented by Panda et al. due to the consideration of an illuminating source composed of both diffuse and oblique collimated flux, which better represents how sunlight interacts with an algal suspension in natural environments. In contrast, Panda et al. used a vertical collimated flux in addition to a diffuse flux, which may not accurately reflect the behaviour of photo-gyro-gravitactic algae in their natural habitat. However, it is important to note that there is currently no quantitative analysis or experimental data available to directly compare with the theoretical developments presented in the proposed model.

	\section*{Data Availability}
The article includes the necessary information and analysis to support the conclusions and results reported.
	\nocite{*}
	\bibliography{anisotropic_rotation}
	
\end{document}